\documentclass{article}

\usepackage{amsthm,amssymb,latexsym,amsfonts,amsmath}
\usepackage{url}
\usepackage{epsfig}
\usepackage{subfigure}
\newtheorem{theorem}{Theorem}[section]
\newtheorem{lemma}[theorem]{Lemma}
\newtheorem{problem}[theorem]{Problem}
\newtheorem{proposition}[theorem]{Proposition}

\newtheorem{definition}[theorem]{Definition}

\newtheorem{remark}[theorem]{Remark}

\newcommand{\R}{{\mathbb{R}}}

\newcommand{\N}{{\mathbb{N}}}

\newcommand{\ie}{{\it i.e.}}

\newcommand{\ea}{{\it et al }}

\begin{document}


\title{An ISS self-triggered implementation\\ of linear controllers.
\thanks{This work has been partially funded by the NSF grants $0841216$ and
 		$0820061$, the Spanish Ministry of Science and Education/UCLA
 		fellowship LA$2004-0003$ and a Mutua Madrile\~na Automovil\'istica
 		scholarship.}}
\author{Manuel Mazo Jr., Adolfo Anta and Paulo
\thanks{M. Mazo Jr, Adolfo Anta and P. Tabuada are with the Department of
		Electrical Engineering, Henry Samueli School of Engineering and Applied
		Sciences, University of California, Los Angeles, CA 90095-1594,
		{\tt\small \{mmazo,adolfo,tabuada\}@ee.ucla.edu}}%
}
\maketitle

%
%

%

\begin{abstract}                          
Nowadays control systems are mostly implemented on digital platforms and,
increasingly, over shared communication networks. Reducing resources
(processor utilization, network bandwidth, etc.) in such implementations
increases the potential to run more applications on the same hardware. We
present a self-triggered implementation of linear controllers that reduces the amount of controller updates necessary to retain
stability of the closed-loop system. Furthermore, we show that the proposed
self-triggered implementation is robust against additive disturbances and
provide explicit guarantees of performance. The proposed technique exhibits
an inherent trade-off between computation and potential savings on actuation.
\end{abstract}


%
%

\section{Introduction}
\label{sec:intro}

The majority of control systems are nowadays implemented on digital
platforms equipped with microprocessors capable of running real-time operating
systems. This creates the possibility of sharing the computational
resources between control and other kinds of applications
thus reducing the deployment costs of complex control systems.
Many control systems are also implemented over shared communication
media making it necessary to share the communication medium.
The concept of self-triggered control was introduced by Velasco
\ea~in~\cite{velasco2003stt} to take advantage of the possibility (or
necessity) of sharing resources. The key idea of self-triggered control is to
compute, based on the current state measurements, the next instant of time at
which the control law is to be recomputed. In between updates of the
controller the control signal is held constant, and the appropriate generation
of the update times guarantees the stability of the closed-loop system. 
Under a periodic implementation, the control law is executed every $T$ units of time, regardless of the current state of the plant. Hence this period $T$ has to be chosen in order to guarantee stability under all possible operating conditions. On the other hand, under self-triggered implementations the time between updates is a function of the state, and thus less control executions are expected. On the other hand, the intervals of time in which no attention is
devoted to the plant pose a new concern regarding the robustness of self-triggered implementations.

The contribution of this paper is to describe a self-triggered implementation
for linear systems, in which the times between controller updates are as large
as possible so as to enforce desired levels of performance subject to the
computational limitations of the digital platform. By increasing the available
computational resources, the performance guarantees improve while the
number of controller executions is reduced. Hence, the proposed technique
reduces the actuation requirements (and communication, in networked systems)
in exchange for computation.
Furthermore, we also show
that the proposed self-triggered implementation results in an exponentially
input-to-state stable closed-loop system with corresponding gains
depending on the available computational resources. A preliminary version of these results appeared
in the conference papers~\cite{MazoJr:2008p1343} and~\cite{MazoJr:CDC2009}.

Several self-triggered implementations
have been proposed in the last years, both for linear~\cite{wang09} and
non-linear~\cite{Anta:2008p2867} plants. The latter when applied to linear
systems degenerates into a periodic implementation, while~\cite{wang09} makes
use of very conservative approximations. In contrast with those two
techniques the approach followed in the present work provides large
inter-execution times for linear systems by not requiring a continuous decay
of the Lyapunov function in use, much in the spirit of~\cite{wang08lt}. 
Computing exactly the maximum allowable inter-execution times guaranteeing
stability requires the solution of transcendental equations for which
closed form expressions do not exist. Our proposal computes approximations
of these maximum allowable inter-execution times while providing stability
guarantees. The idea advocated in this paper, trading
communication/actuation for computation, was already explored
in~\cite{Yook:2002}. However, their approach is aimed at loosely coupled
distributed systems, where local actuation takes place continuously and
communication between subsystems is reduced by means of state estimators.
In the analysis of robustness of the proposed implementation the authors were
influenced by the approach followed in~~\cite{Nesic:2004p1869} and~\cite{Kellett:2004p168}.
Finally, the notion of input-to-state stability~\cite{Sontag:2008p2453} is
fundamental in the approach followed in the present paper.

\section{Notation}
\label{sec:prelim}

We denote  
by $\R^+$ the positive
real numbers. We also use
$\R_0^+=\R^+\cup\lbrace0\rbrace$. The usual Euclidean ($l_2$) vector norm is
represented by $|\cdot|$. When applied to a matrix $|\cdot|$ denotes the
$l_2$ induced matrix norm. A matrix
$P\in\R^{m\times m}$ is said to be positive definite, denoted $P>0$, whenever
$x^TPx>0$ for all $x\neq 0$, $x\in\R^m$, and a matrix $A$ is said to be Hurwitz
when all its eigenvalues have strictly negative real part. We denote by $I$ the
identity matrix. By $\lambda_m(P),\lambda_M(P)$ we denote the minimum and
maximum eigenvalues of $P$ respectively. A function $\gamma:[0,\infty[\to\R^+_0$,
is of class~$\mathcal{K}_\infty$ if it is continuous, strictly
increasing, $\gamma(0)=0$ and $\gamma(s)\to\infty$ as $s\to\infty$. Given an
essentially bounded function $\delta:\R^+_0\to\R^m$ we denote by $\Vert \delta
\Vert_\infty$ its $\mathcal{L}_\infty$ norm, \ie, $\Vert \delta
\Vert_\infty=(ess)\sup_{t\in\R^+_0}\lbrace{|\delta(t)|\rbrace}<\infty$. We
consider linear systems described by the differential equation
$\frac{d}{dt}\xi=A\xi+B\chi+\delta $ with inputs $\chi:\R^+_0\to\R^l$ and
$\delta:\R^+_0\to\R^p$ essentially bounded piecewise continuous functions of
time. The input $\chi$ will be used to denote controlled inputs, while
$\delta$ will denote disturbances. We refer to such systems as {\it control
systems}. Solutions of a control system with initial condition $x$ and inputs
$\chi$ and $\delta$, denoted by $\xi_{x\chi\delta}$, satisfy:
$\xi_{x\chi\delta}(0)=x$ and
$\frac{d}{dt}\xi_{x\chi\delta}(t)=A\xi_{x\chi\delta}(t)+B\chi(t)+\delta(t)$
for almost all $t\in\R^+_0$. The notation will be relaxed by dropping the
subindex when it does not contribute to the clarity of exposition. A linear
feedback law for a linear control system is a map $u=Kx$; we will sometimes
refer to such a law as a {\it controller} for the system.

\begin{definition}[Lyapunov function]
A function \mbox{$V:\R^m\to\R^+_0$}, is said to
be a Lyapunov function for a linear system $\dot{\xi}=A\xi$ if $\forall\,x\in\R^m$
$\underline{\alpha}(|x|)\leq V(x)\leq\overline{\alpha}(|x|)$ for some
$\underline{\alpha},\overline{\alpha}\in\mathcal{K}_\infty$,
and there exists $\lambda\in\R^+$ such that for every
$x\in\R^m$:
\begin{equation}\notag
\frac{\partial V}{\partial x}Ax\leq-\lambda V(x).
\end{equation}
\end{definition}
We will refer to $\lambda$ as the {\it rate of decay} of the Lyapunov function.
In what follows we will consider functions of the form
$V(x)=(x^TPx)^\frac{1}{2}$,
in which case $V$ is a Lyapunov function for system
\mbox{$\dot{\xi}=A\xi$} if and only if $P>0$ and
$A^TP+PA\leq-2\lambda I$ for some $\lambda\in\R^+$, the rate of decay.

\begin{definition}[EISS]
A control system $\dot{\xi}=A\xi+\delta$ is said to be
{\it exponentially input-to-state stable} (EISS) if there exists
$\lambda\in\R^+$, $\sigma\in\R^+$ and $\gamma\in\mathcal{K}_\infty$ such
that for any $t\in\R^+_0$ and for all $x\in\R^m$:
\begin{equation}\notag
|\xi_{x\delta}(t)|\leq\sigma|x|e^{-\lambda t}+\gamma(\Vert
\delta \Vert_\infty).
\end{equation}
\end{definition}
We shall refer to $(\beta,\gamma)$, where $\beta(s,t)=s\sigma e^{-\lambda t}$,
as the {\it EISS gains} of the EISS estimate.
If no disturbance is present, \ie, $\delta=0$, an EISS
system is said to be {\it globally exponentially stable} (GES).

\section{A self-triggered implementation for stabilizing linear controllers.}
\label{sec:proposal}

Consider the sampled-data system: \begin{eqnarray}
\label{eq:HybSys1}
\dot{\xi}(t)&=&A\xi(t)+B\chi(t)+\delta(t)\\
\label{eq:HybSys1_contr}
\chi(t)&=&K\xi(t_k), \; t\in[t_k,t_{k+1}[
\end{eqnarray}
where $\lbrace{t_k\rbrace}_{k\in\N}$ is a divergent sequence of update times
for the controller, and $A+BK$ is Hurwitz. The signal
$\delta$ can be used to describe measurement disturbances, actuation
disturbances, unmodeled dynamics, or other sources of uncertainty as
described in~\cite{Kellett:2004p168}.

A self-triggered implementation of the linear stabilizing
controller~(\ref{eq:HybSys1_contr}) for the plant~(\ref{eq:HybSys1}) is given
by a map \mbox{$\Gamma:\R^m\to \R^+$} determining the controller update time $t_{k+1}$
as a function of the state $\xi(t_k)$ at the time $t_k$, \ie,
$t_{k+1}=t_k+\Gamma(\xi(t_k))$. If we denote by $\tau_k$ the inter-execution
times $\tau_k=t_{k+1}-t_k$, we have $\tau_k=\Gamma(\xi(t_k))$. Once the map
$\Gamma$ is defined, the expression {\it closed-loop system} refers to the
sampled-data system~(\ref{eq:HybSys1}) and~(\ref{eq:HybSys1_contr}) with the
update times $t_k$ defined by $t_{k+1}=t_k+\Gamma(\xi(t_k))$.

The problem we solve in this
paper is the following:
\begin{problem}
Given a linear system~(\ref{eq:HybSys1}) and a linear stabilizing
controller~(\ref{eq:HybSys1_contr}), construct a self-triggered
implementation $\Gamma:\R^m\to\R^+$ of~(\ref{eq:HybSys1_contr}) that renders
EISS the closed-loop system defined
by~(\ref{eq:HybSys1}),~(\ref{eq:HybSys1_contr}), while enlarging the inter-execution times.
\end{problem}
In order to formally define the self-triggered implementation proposed in this paper, we
need to introduce two maps:
\begin{itemize}
  \item $h_c$, a continuous-time output map and
  \item $h_d$, a discrete-time version of $h_c$.
\end{itemize}

Let $V$ be a Lyapunov function of the form  \linebreak$V(x)=(x^TPx)^\frac{1}{2}$ for
\mbox{$\dot{\xi}=(A+BK)\xi$}, with rate of decay $\lambda_o$. The output
map $h_c:\R^m\times\R_0^+\to\R_0^+$ is defined by:
\begin{equation}
\label{eq:et_func}
h_c(x,t):=V(\xi_x(t))-V(x)e^{-\lambda t}
\end{equation}
for some $0<\lambda<\lambda_o$. Note that by enforcing:
\begin{equation}
\label{DesiredIneq}
h_c(\xi_x(t_k),t)\leq0,\quad \forall t\in[0,\tau_k[\,\,\forall k\in\N 
\end{equation}
the closed-loop system satisfies:
\begin{equation*}
V(\xi_{x}(t))\leq V(x)e^{-\lambda t}, \quad\forall t\in\R^+_0\,\, \forall
x\in\R^m
\end{equation*}
which implies exponential stability of the closed-loop system in the absence of
disturbances, \ie, when $\delta(t)=0$ for all $t\in\R^+_0$. 

Our objective is to construct a self-triggered implementation
enforcing~(\ref{DesiredIneq}). Since no digital implementation can
check~(\ref{DesiredIneq}) for all $t\in[t_k,t_{k+1}[$, we consider instead the
following discrete-time version of~(\ref{DesiredIneq}) based on a sampling
time $\Delta\in \R^+$: 
$$h_d(\xi_x(t_k),n):=h_c(\xi_x(t_k),n\Delta)\leq0\quad\forall
n\in\left[0,\left\lceil\frac{\tau_k}{\Delta}\right\rceil\right[,$$
and for all $k\in \N$. This results
in the following self-triggered implementation, first introduced by the
authors in~\cite{MazoJr:2008p1343}, where we use
$N_{\min}:=\left\lfloor{\tau_{\min}}/{\Delta}\right\rfloor$, \mbox{$
N_{\max}:=\left\lfloor{\tau_{\max}}/{\Delta}\right\rfloor $}, and $\tau_{\min}$
and $\tau_{\max}$ are design parameters. A similar approach was followed in~\cite{heemels08} in the context of event-triggered control.
\begin{definition}
\label{def:st}
The map $\Gamma_d:\R^n\to\R^+$ defined by: 
\begin{eqnarray}\notag
\Gamma_d(x)&:=&\max\lbrace{\tau_{\min},
n_k\Delta\rbrace}\\\notag
n_k&:=&\max_{n\in\N}\lbrace{n\leq N_{\max}
\vert h_d(x,s)\le 0, s=0,\hdots,n\rbrace}
\end{eqnarray}
prescribes a self-triggered implementation of the linear stabilizing
controller~(\ref{eq:HybSys1_contr}) for plant~(\ref{eq:HybSys1}).
\end{definition}

Note that the role of $\tau_{\min}$ and $\tau_{\max}$ is to enforce explicit
lower and upper bounds, respectively, for the inter-execution times of the
controller. The upper bound enforces robustness of the 
implementation and limits the computational complexity.

\begin{remark}
Linearity of~(\ref{eq:HybSys1}) and~(\ref{eq:HybSys1_contr}) enables us to
compute $h_d$ as a quadratic function of $\xi(t_k)$. Moreover, through a
Veronese embedding we can implement the self-triggered policy described in
Definition~\ref{def:st} so that its computation has space complexity
\mbox{$q\frac{m(m+1)}{2}$} and time complexity
\mbox{$q+(2q+1)\frac{m(m+1)}{2}$} where \mbox{$q:=N_{\max}-N_{\min}$}. For
reasons of space we omit these details. They can be
found in~\cite{MazoJr:2008p1343}.
\end{remark}

\section{Main results}
\label{sec:analysis}
The proofs of all the results reported in this section can be found in the
Appendix. The following functions will be used to define EISS-gains: 
$$
\rho_P:=\left(\frac{\lambda_M(P)}{\lambda_m(P)}\right)^\frac{1}{2},\;
\gamma_{P,T}(s):=s\frac{\lambda_M(P)}{\lambda_m^\frac{1}{2}(P)}\int_0^{T}\!\!\!
\vert e^{Ar}\vert dr. 
$$

We start by establishing a result explaining how the design parameter
$\tau_{\min}$ should be chosen. The function $\Gamma_d$ can be seen as a
discrete-time version of the function \mbox{$\Gamma_c:\R^m\to\R_0^+$} defined by:
\begin{equation}
\label{eq:intersamples_c}
\Gamma_c(x):=\max_{\tau\in\R^+_0}\lbrace{\tau\leq\tau_{\max} |
h_c(x,s)\leq0,\forall s\in[0,\tau]\rbrace}.
\end{equation}
If we use $\Gamma_c$ to define an ideal self-triggered implementation, the
resulting inter-execution times are no smaller than $\tau^*_{\min}$ which can
be computed as detailed in the next result.

\begin{lemma}
\label{lmm:min_time}
The inter-execution times generated by the self-triggered implementation
in~(\ref{eq:intersamples_c}) are lower bounded by:
\begin{equation}
\label{eq:min_time}
\tau^*_{\min} = \min \{\tau \in \R^+ : \det M(\tau) = 0\}
\end{equation}
where:
\begin{eqnarray*}
M(\tau) &: =& C(e^{F^T\tau}C^TPCe^{F\tau} - C^TPCe^{-\lambda \tau})C^T,\\
F&:=& \left[ \begin{array}{ccc}  \phantom{-}A+BK &
\phantom{-}BK \\ -A-BK & -BK \end{array} \right],\; C:=[I\; 0].
\end{eqnarray*}
\end{lemma}
The computation of $\tau^*_{\min}$ described in Lemma~\ref{lmm:min_time} can
be regarded as a formal procedure to find a sampling period for periodic
implementations (also known as maximum allowable time interval or MATI). It should be contrasted with the frequently used ad-hoc rules of thumb~\cite{astrom90}, \cite{handbook-networked-control} (which do not provide stability guarantees).
Moreover, an analysis similar to the one in the proof of this lemma can also be
applied, $\mathit{mutatis}$ $\mathit{mutandis}$, to other Lyapunov-based
triggering conditions, like the ones appearing in~\cite{tabuada07}
and~\cite{wang09}. Notice that the self-triggered approach always renders times no smaller than the periodic implementation, since under a periodic implementation the controller needs to be executed every $\tau^*_{\min}$ (in order to guarantee performance under all possible operating points).

The second and main result establishes EISS of the proposed self-triggered
implementation and details how the design parameters $\tau_{\min},
\tau_{\max}, \Delta,$ and $\lambda$ affect the EISS-gains.

\begin{theorem}
\label{thm:ISS}
If $\tau_{\min}\leq\tau^*_{\min}$, the self-triggered implementation in
Definition~\ref{def:st} renders the closed-loop system
EISS with gains $(\beta,\gamma)$ given by:
\begin{small}
\begin{eqnarray*}
\beta(s,t)&:=&\rho_P g(\Delta,N_{\max})e^{-\lambda t}s,\\
\gamma(s)&:=&\gamma_{P,N_{\max}\Delta}(s)\frac{\lambda_m^{-\frac{1}{2}}(P)
g(\Delta,N_{\max})}{1-e^{-\lambda \tau_{\min}}}+\gamma_{I,N_{\max}\Delta}(s)
\end{eqnarray*}
where:
\begin{eqnarray*}
g(\Delta,N_{\max}):=\rho_P\left(e^\frac{(\rho+2\lambda)\mu\Delta}{\mu-\rho}+\right.\qquad\qquad\qquad\\
\left.+e^{2\lambda(N_{\max}-1)\Delta}\left(
e^\frac{(\rho+2\lambda)\mu\Delta}{\mu-\rho}-e^\frac{2\lambda\mu\Delta}{\mu-\rho}\right)\right)^\frac{1}{2},
\end{eqnarray*}
\vspace{-0.7cm}
\begin{eqnarray*}
\rho&:=&\lambda_M(G),\; \mu:=\lambda_m(G),\\ 
G&:=&\left[\begin{array}{cc}
P^{\frac{1}{2}}AP^{-\frac{1}{2}}+(P^{\frac{1}{2}}AP^{-\frac{1}{2}})^T\;\; &
P^{\frac{1}{2}}BKP^{-\frac{1}{2}} \\ (P^{\frac{1}{2}}BKP^{-\frac{1}{2}})^T & 0
\end{array}\right].
\end{eqnarray*}
\end{small}
\end{theorem}

Note that while $\tau_{\min}$ is constrained by $\tau_{\min}^*$, $\tau_{\max}$
can be freely chosen. However, by enlarging $\tau_{\max}$ (and thus $N_{\max}$) we
are degrading the EISS-gains. It is also worth noting that by enlarging
$\tau_{\max}$ one can allow longer inter-execution times, and compensate the
performance loss by decreasing $\Delta$, at the cost of performing more computations.

Let us define the maximum exact inter-execution
time from $x$ as 
$\tau^*(x):=\min\lbrace \Gamma_c(x), \tau_{\max}\rbrace$, where the upper bound is required to obtain robustness against disturbances.
The third and final result states that the proposed self-triggered
implementation is optimal in the sense that it generates the longest possible
inter-execution times given enough computational resources. Hence, by enlarging
the inter-execution times we are effectively trading actuation for computation.
The proof of the following proposition follows from the proof of
Theorem~\ref{thm:ISS}.

\begin{proposition}
The inter-execution times provided by the self-triggered implementation in
Definition~\ref{def:st} are bounded from below as follows:
$$\Gamma_d(x)\geq\tau^*(x)-\Delta,\;\forall\, x\in\R^m.$$
\end{proposition}
Note that even if $\Gamma_d(x)\geq\tau^*(x)$ the performance guarantees
provided in Theorem~\ref{thm:ISS} still hold.

\begin{remark}
When implementing self-triggered policies on digital platforms several
issues related to real-time scheduling need to be addressed. For a discussion
of some of these issues we refer the readers to~\cite{Anta:2009}. Here,
we describe the minimal computational requirements for the proposed
self-triggered implementation under the absence of other tasks. There are
three main sources of delays: measurement, computation, and actuation. Since
the computation delays dominate the measurement and actuation delays, we focus
on the former. The computation of~$\,\Gamma_d$ is divided in two steps: a
preprocessing step performed once by execution, and a running step performed
$n$ times when computing $h_d(x,n)$. The preprocessing step has time
complexity $(m^2+m)/2$ and the running step has time complexity $m^2+m$. If we
denote by $\tau_c$ the time it takes to execute an instruction in a given
digital platform, the self-triggered implementation can be executed if: 
$$
\frac{3}{2}(m^2+m)\tau_c\leq\tau_{\min},\;(m^2+m)\tau_c\leq\Delta.
$$
The first inequality imposes a minimum processing speed for the digital
platform while the second equality establishes a lower bound for the choice of
$\Delta$.
\end{remark}

We refer the interested reader to~\cite{MazoJr:2008p1343}
and~\cite{MazoJr:CDC2009} for numerical examples illustrating the proposed
technique and the guarantees it provides.

\section{Conclusions}

This paper described a self-triggered implementation of stabilizing feedback
control laws for linear systems. The proposed technique guarantees exponential
input-to-state stability of the closed-loop system with respect to additive
disturbances. Furthermore, the proposed self-triggered implementation allows
the tuning of the resulting performance and complexity through the selection
of the parameters $\Delta$, $\lambda$, $\tau_{\min}$ and $\tau_{\max}$. The
performance guarantees can be improved and the inter-execution times
enlarged by increasing the computational complexity of the implementation.

\bibliographystyle{IEEEtran}
\bibliography{aut_self_iss_review_final} 

\section{Appendix: Proofs}

\begin{proof}[Proof of Lemma~\ref{lmm:min_time}]
It can be verified that $h_c$ satisfies
$h_c(x,0)=0$ and $\frac{\partial}{\partial
t}\Big\vert_{t=0}h_c(x,t)<0$,$\forall\,x\in\R^m$,
which, by continuity of $h_c$, implies the existence of some \mbox{$\tau^*_{\min}(x)>0$}
such that $\Gamma_c(x)\geq \tau^*_{\min}(x)$.
Let us define the variables $\eta(t)=\xi(t)-\xi(t_k)$,
$t\in[t_k,t_{k+1}[$ and \mbox{$\zeta=[\xi^T\; \eta^T]^T$}. Note that at the
controller update times $\eta(t_k)=0$. Under this new notation,
system~(\ref{eq:HybSys1}) with controller~(\ref{eq:HybSys1_contr}), in the absence of disturbances, can be
rewritten as $\dot{\zeta}(t) = F\zeta(t)$ with solution $\zeta_{y}(t) = e^{Ft}
y$, $y=[x^T\; 0^T]^T$. Let us denote by $\hat{h}_c$ the map
$\hat{h}_c(y,t)=V(C\zeta_y(t))-V(Cy)e^{-\lambda t}$.

While it is not possible to find $\Gamma_c$ in closed form, we can find its
minimum value by means of the Implicit Function Theorem. Differentiating
$\phi(x)=\hat{h}_c(y,\Gamma_c(x))=0$ with respect to the initial condition $x$ we obtain:
\begin{eqnarray*}
\frac{d\phi}{dx} &=& \frac{\partial \hat{h}_c}{\partial
t}\Big\vert_{t=\Gamma_c(x)}\frac{d\Gamma_c}{dx}+\frac{\partial
\hat{h}_c}{\partial y}\frac{d y}{dx}=0.
\end{eqnarray*}
The extrema of the map $\Gamma_c$ are defined by the following equation:
\begin{equation*}
\frac{d\Gamma_c}{dx}=-\left( \frac{\partial \hat{h}_c}{\partial
t}\Big\vert_{t=\Gamma_c(x)} \right)^{-1}\left(\frac{\partial \hat{h}_c}{\partial y}\frac{d
y}{dx}\right)=0.
\end{equation*}
Hence, the extrema of $\Gamma_c$
satisfy either \mbox{$\frac{\partial \hat{h}_c}{\partial y}\frac{d
y}{dx}(x,t)=0$} \mbox{for some $t\in\R^+$}
or \mbox{$\frac{d\Gamma_c}{dx}(x)=0\; \wedge$} \mbox{$\frac{\partial
\hat{h}_c}{\partial t}\Big\vert_{t=\Gamma_c(x)}=0$}. The latter case corresponds to situations in
which for some $x$ the map $\hat{h}_c$ reaches zero exactly at an extremum, and
thus can be disregarded as violations of the condition $h_c(t_k,t) \leq 0$.
Combining $\frac{\partial \hat{h}_c}{\partial y}\frac{d
y}{dx}(\tau,x)=0$ into matrix form we obtain:
\begin{equation}
\label{final_cond}
M(\tau) x = 0.
\end{equation}
The solution to this equation provides all extrema of the map
$\Gamma_c(x)$ that incur a violation of $h_c(x,t) \leq 0$. Thus, the minimum
$\tau$ satisfying~(\ref{final_cond}) corresponds to the smallest
time at which \mbox{$h_c(x,\tau)=0$}, \mbox{$\frac{\partial}{\partial
t}h_c(x,t)\Big\vert_{t=\tau}>0$} can occur. Since the left hand side
of~(\ref{final_cond}) is linear in $x$, it is sufficient to check when the
matrix has a nontrivial nullspace. Hence the equality~(\ref{eq:min_time}).
\end{proof}

We introduce now a Lemma that will be used in the proof of
Theorem~\ref{thm:ISS}.

\begin{lemma}
\label{lemma:disturb}
Consider system~(\ref{eq:HybSys1}) and a positive definite function
$V(x)=\left(x^TPx\right)^\frac{1}{2},\; P>0$. For any given \mbox{$0\leq
T<\infty$} the following bound holds:
\begin{equation*}
V(\xi_{x\chi\delta}(t))\leq V(\xi_{x\chi 0}(t))+\gamma_{P,T}(\Vert \delta
\Vert_\infty), \forall t\in[0,T].
\end{equation*}
\end{lemma}
\begin{proof}
Applying the triangular inequality and using Lipschitz continuity of $V$ we
have: 
\begin{eqnarray}\notag
V(\xi_{x\chi\delta}(t))&=&|V(\xi_{x\chi 0}(t))+V(\xi_{x\chi\delta}(t))-V(\xi_{x\chi 0}(t))|\\\notag
&\leq& V(\xi_{x\chi
0}(t))+\frac{\lambda_M(P)}{\lambda_m^\frac{1}{2}(P)}|\xi_{x\chi\delta}(t)-\xi_{x\chi
0}(t)|.
\end{eqnarray}
Integrating the dynamics of $\xi$ and after
applying H{\"o}lder's inequality one can conclude that:
\begin{equation*}
|\xi_{x\chi\delta}(t)-\xi_{x\chi 0}(t)|\leq \int_0^t \vert e^{Ar}\vert
dr\Vert \delta \Vert_\infty.
\end{equation*}
And thus for all $t\in[0,T]$:
\begin{eqnarray*}
V(\xi_{x\chi\delta}(t))&\leq& 
V(\xi_{x\chi 0}(t))+\frac{\lambda_M(P)}{\lambda_m^\frac{1}{2}(P)}\int_0^T \vert
e^{Ar}\vert dr\Vert \delta \Vert_\infty.
\end{eqnarray*}
\end{proof}

\begin{proof}
[Proof of Theorem~\ref{thm:ISS}]

We start by proving that in the absence of disturbances the following bound
holds:
\begin{equation}
\label{eq:bound_nd}
|\xi_x(t_k+\tau)|\leq
g(\Delta,N_{\max})|\xi_x(t_k)|e^{-\lambda\tau},\;\forall\,\tau\geq 0.
\end{equation}
Let $W(x)=x^TPx$ and use $W(t)$ to denote
$W(\xi_{x_k}(t))$, with $\xi$ determined
by~(\ref{eq:HybSys1}),~(\ref{eq:HybSys1_contr}), and
$\tau_k=\Gamma_d(\xi(t_k))$. By explicitly computing $\dot{W}(t)$ one
obtains:
\begin{eqnarray*}
&&\dot{W}(t)=\\
&&\left[(P^\frac{1}{2}\xi(t))^T\;(P^\frac{1}{2}\xi(t_k))^T\right]G\left[(P^\frac{1}{2}\xi(t))^T\;(P^\frac{1}{2}\xi(t_k))^T\right]^T,
\end{eqnarray*}
for $t\in[t_k,t_{k+1}[$, and thus the following bounds hold:
\begin{equation*}
\mu\left(W(t)+W(t_k)\right)\leq \dot{W}(t) \leq\rho\left(W(t)+W(t_k)\right).
\end{equation*}
for $t\in[t_k,t_{k+1}[$.
After integration, one can bound the trajectories of $W(t)$, when $t+s$ belongs
to the interval $[t_k,t_{k+1}[$, as:
\begin{eqnarray*}
W(t+s)&\leq & e^{\rho s}W(t) + W(t_k)\left(e^{\rho s} - 1\right),\\
W(t+s)&\geq & e^{\mu s}W(t) + W(t_k)\left(e^{\mu s} - 1\right).
\end{eqnarray*}
Let us denote $t_k+n\Delta$ by $r_n$ for succinctness of the expressions that
follow. An upper bound for $W(t)$ valid for $t\in[r_n, r_{n+1}[$
is then provided by:
\vspace{-0.7cm}
\begin{small}
\begin{eqnarray*}
&&W(r_n+s)\leq\\
&&\left\lbrace\begin{array}{ll}
e^{\rho s}W(r_n)+W(t_k)(e^{\rho s}-1),&s\in[0,s^*]\\
e^{\mu(s-\Delta)}W(r_n+\Delta)+W(t_k)(e^{\mu(s-\Delta)}-1),&s\in[s^*,\Delta[.
\end{array}\right.
\end{eqnarray*}
\end{small}
\vspace{-0.7cm}

The maximum for the bound of
$W(r_n+s)$ for $s\in[0,\Delta[$, is attained at the point at
which the two branches of the bound meet, \ie at $s=s^*$, as the first branch
is monotonically increasing in $s$, and the second branch monotonically
decreasing. The point $s^*$ can be computed as:
\begin{equation*}
s^*=\frac{1}{\rho-\mu}\log\left(\frac{W(r_{n+1})+W(t_k)}{W(r_n)+W(t_k)}\right)+\frac{\mu\Delta}{\mu-\rho}
\end{equation*}
and thus $W(r_n+s^*)$ can be bounded as:
\vspace{-0.7cm}
\begin{small}
\begin{eqnarray*}
&&W(r_n+s^*)\leq-W(t_k)+\\
&&\quad\;\;
e^\frac{\rho\mu\Delta}{\mu-\rho}\left(\left(W(r_n)+W(t_k)\right)^\frac{\mu}{\mu-\rho}
\left(W(r_{n+1})+W(t_k)\right)^\frac{-\rho}{\mu-\rho}\right)
\end{eqnarray*}
\end{small}
\vspace{-0.7cm}

which is monotonically increasing on $W(r_n)$, $W(r_{n+1})$, 
and $W(t_k)$. If $S(t)=W(t_k)e^{-2\lambda(t-t_k)}$, it follows:
\vspace{-0.7cm}
\begin{small}
\begin{eqnarray*}
&&W(r_n+s^*)\leq -S(t_k)+\\
&&\qquad
e^\frac{\rho\mu\Delta}{\mu-\rho}\left(\left(S(r_n)+S(t_k)\right)^\frac{\mu}{\mu-\rho}
\left(S(r_{n+1})+S(t_k)\right)^\frac{-\rho}{\mu-\rho}\right)
\end{eqnarray*}
\end{small}
\vspace{-0.7cm}

where we used the fact that, if $\tau_{\min}\leq\tau^*_{\min}$, $\Gamma_d$
enforces (in the absence of disturbances) $W(r_n)\leq S(r_n)$
for all $n\in\N$, $n\leq n_k$. 
From the previous expression we can obtain
\mbox{$W(r_n+s^*)\leq\tilde{g}(\Delta,n)S(r_n+s^*)$} where:
\begin{scriptsize}
\begin{eqnarray}\notag
&&\tilde{g}(\Delta,n)= - e^{2\lambda(n\Delta+s^*)} +\\\notag
&&e^\frac{\rho\mu\Delta}{\mu-\rho}\left(e^{2\lambda
s^*}(1+e^{2\lambda n\Delta}) \right)^\frac{\mu}{\mu-\rho}\!\!
\left(e^{-2\lambda(\Delta-s^*)}+
e^{2\lambda(n\Delta+s^*)}\right)^\frac{-\rho}{\mu-\rho}\!\!.
\end{eqnarray}
\end{scriptsize}

The value of $s^*$ can be further
bounded to obtain a simpler expression:
\begin{equation*}
s^*\leq\frac{\mu\Delta}{\mu-\rho}.
\end{equation*}
Using this bound for $s^*$ and letting $n$ take its maximum possible value
$n=N_{\max}-1$, the following chain of inequalities holds:
$$\rho_{P}\tilde{g}(\Delta,n)^{\frac{1}{2}}\leq\rho_{P}\tilde{g}(\Delta,N_{\max}-1)^{\frac{1}{2}}\leq
g(\Delta,N_{\max})$$ 
for all $n\in[0,N_{\max}]$, which leads to the bound:
\begin{equation}
\label{eq:bndg}
W^{\frac{1}{2}}(t)\leq \rho_{P}^{-1}g(\Delta,
N_{\max})S^{\frac{1}{2}}(t).
\end{equation} 
Note that~(\ref{eq:bndg}) does not depend on $t_k$ or $n$.
Finally, apply the bounds:
\begin{equation}
\label{eq:boundsV}
\lambda_m^{\frac{1}{2}}(P)|x|\leq
V(x)=\sqrt{x^TPx}\leq \lambda_M^{\frac{1}{2}}(P)|x|.
\end{equation}
to obtain~(\ref{eq:bound_nd}).
From Lemma~\ref{lemma:disturb}, and the condition enforced by the
self-triggered implementation we have:
\begin{equation*}
V(\xi(t_{k+1}))\leq V(\xi(t_k))e^{-\lambda \tau_k}+\gamma_P(\Vert\delta
\Vert_\infty).
\end{equation*}
Iterating the previous equation it follows:
\begin{eqnarray*}
V(\xi(t_k))\!&\leq&\!e^{-\lambda(t_k-t_o)}
V(\xi(t_o))+\gamma_P(\Vert\delta \Vert_\infty)\sum_{i=0}^{k-1}e^{-\lambda t_{\min}i}\\
\!&\leq&\!e^{-\lambda(t_k-t_o)}V(\xi(t_o))+\gamma_P(\Vert\delta
\Vert_\infty)\sum_{i=0}^{\infty}e^{-\lambda t_{\min}i}\\
\!&=&\!e^{-\lambda(t_k-t_o)}V(\xi(t_o))+\gamma_P(\Vert\delta
\Vert_\infty)\frac{1}{1-e^{-\lambda t_{\min}}}.
\end{eqnarray*}
Assuming, without loss of generality, that $t_o=0$, the following bound
also holds:
\begin{equation}
\label{eq:proofISS1}
|\xi_x(t_k)|\leq \rho_P|x|e^{-\lambda
t_k}+\lambda_m^{-\frac{1}{2}}(P)\frac{\gamma_P(\Vert\delta
\Vert_\infty)}{1-e^{-\lambda t_{\min}}}
\end{equation}
where we used~(\ref{eq:boundsV}). 
From~(\ref{eq:bound_nd}) and
Lemma~\ref{lemma:disturb} one obtains:
\begin{equation}
\label{eq:proofISS2}
|\xi_x(t_k+\tau)|\leq
g(\Delta,N_{\max})|\xi_x(t_k)|e^{-\lambda\tau}+\gamma_I(\Vert\delta\Vert_\infty),
\end{equation}
for all $\tau\in[0,N_{\max}\Delta]$.
Combining~(\ref{eq:proofISS1})~and~(\ref{eq:proofISS2}) results in:
\begin{eqnarray*}
|\xi_x(t_k+\tau)|&\leq& g(\Delta,N_{\max})\rho_P|x|e^{-\lambda
(t_k+\tau)}\\
&+&e^{-\lambda\tau}\gamma_P(\Vert\delta\Vert_\infty)\frac{\lambda_m^{-\frac{1}{2}}(P)
g(\Delta,N_{\max})}{1-e^{-\lambda t_{\min}}}\\
&+&\gamma_I(\Vert\delta\Vert_\infty),
\end{eqnarray*}
and after denoting $t_k+\tau$ by $t$ we can further bound:
\begin{eqnarray*}
|\xi_x(t)|&\leq& g(\Delta,N_{\max})\rho_P|x|e^{-\lambda t}\\
&+&\gamma_P(\Vert\delta\Vert_\infty)\frac{\lambda_m^{-\frac{1}{2}}(P)
g(\Delta,N_{\max})}{1-e^{-\lambda t_{\min}}}+\gamma_I(\Vert\delta\Vert_\infty),
\end{eqnarray*}
which is independent of $k$ and concludes the proof.
\end{proof}
\end{document}